\documentclass[12pt,leqno]{article}

\usepackage{amsmath}
\usepackage{amsfonts}
\usepackage{amssymb}
\usepackage{amsthm}
\usepackage{latexsym}
 \sloppypar \textheight 21cm
\newtheorem{theorem}{Theorem}[section]

\newtheorem{lemma}[theorem]{Lemma}
\newtheorem{remark}[theorem]{Remark}
\newtheorem{proposition}[theorem]{Proposition}
\newtheorem{cor}[theorem]{Corollary}
\newtheorem{example}[theorem]{Example}
\newcommand{\supp}{\text{supp}}

\def\N{\mathbb N}

\def\DE{\Delta}

\def\va{\varepsilon}

\def\sp#1#2{\langle{#1},{#2}\rangle}

\def\N{\mathbb{N}}

\def\N{{\mathbb N}}

\def\L{{\cal L}}

\def\val#1{\vert #1\vert}

\def\no#1{\Vert #1\Vert }

\def\inv{^{-1}}

\def\val#1{\vert #1\vert}

\def\sp#1#2{\langle #1,#2\rangle }

\def\L1#1{L^1(#1)}

\def\lef({\left(}
\def\rig){\right)}

\begin{document}

\title{On the connection between sets of operator synthesis and sets
  of spectral
  synthesis for locally compact groups}
\author{Jean Ludwig \vspace{0.1cm}\\
{\footnotesize \sl Department of Mathematics, University of Metz,
}\\
{\footnotesize \sl
 F-57045 Metz, France} \vspace{0.3cm}\\
Lyudmila  Turowska \vspace{0.1cm}\\
{\footnotesize \sl Department of Mathematics, Chalmers University of 
Technology, }\\ 
{\footnotesize \sl SE-412 96 G\"oteborg, Sweden} }
 
\date{}

  \maketitle{}

\begin{abstract} We extend the results by Froelich and Spronk \&
  Turowska on the
connection between operator synthesis and spectral synthesis for
$A(G)$ to  second countable locally compact groups $G$. This gives
us another proof that one-point subset of $G$ is a set of spectral
synthesis and that any closed subgroup is  a set of local spectral
synthesis. Furthermore we show that ``non-triangular'' sets are
strong operator Ditkin sets and we establish  a connection between
operator Ditkin sets and Ditkin sets. These results are applied to
prove that any closed subgroup of $G$  is a local Ditkin set.
\end{abstract}


\section{Introduction}
In \cite{arv} Arveson discovered a connection between the
invariant subspace theory and spectral synthesis. He defined
(operator) synthesis for subspace lattices and proved the failure
of operator synthesis by using the famous  example of Schwartz on
non-synthesizability of the two-sphere $S^2$ for $A({\mathbb
R}^3)$. In \cite{froelich} Froelich made this connection more
precise for separable abelian group. For $G$ a separable compact
group  this relation was obtained in \cite[Theorem~4.6.]{st}. We
generalize these results to second countable locally compact
groups (Theorems~\ref{synthesis} and \ref{ops}). We use the
definition of sets of operator synthesis as defined in \cite{sht}.
We prove that a closed subset $E\subset G$ is set of local
spectral synthesis for $A(G)$ if and only if the diagonal set
$E^*=\{(s,t)\in G\times G\mid st^{-1}\in E\}$ is a set of operator
synthesis with respect to Haar measure. We remark that if $A(G)$
has an (unbounded) approximate identity then any set of local
spectral synthesis is a set of synthesis for $A(G)$ and any
compact set of local spectral synthesis is globally spectral for
any group $G$. We give a new proof that a one-point set is
spectral and any closed subgroup of second countable group is a
set of local spectral synthesis. Using operator synthesizability
of sets of finite width we obtain certain example of sets of
spectral synthesis.

Operator-Ditkin sets have been defined in \cite{sht}. In Section~5
we show first that ``non-triangular'' type sets are strong
operator-Ditkin. The proof is inspired by \cite{drury}. For $G$
second countable group we prove that the diagonal set $E^*\subset
G\times G$ is a strong operator-Ditkin set with respect to Haar
measure then $E\subset G$ is a local Ditkin set for $A(G)$ and
conversely, for any strong Ditkin set $E\subset G$, the set $E^*$ is
operator-Ditkin with respect to Haar measure. As an application we
obtain that any closed subgroup of a second countable group is a
local Ditkin set. This result was known for neutral subgroups of
arbitrary locally compact groups $G$ (\cite{derighetti}) and
amenable groups $G$ (\cite{fkls}).

\section{Preliminaries and notations}
Let $G$ be a locally compact $\sigma$-compact separable group with
 left
 Haar measure $m=dg$. Let $L^p(G)$, $p=1,2$, denote the space of $p$-integrable
functions with norm $||\cdot||_p$ and let $C_c(G)$ denote the
algebra of continuous compactly supported complex-valued functions
on $G$. The convolution algebra $L^1(G)$ is an involutive algebra
with involution defined by $f^*(s)=\Delta^{-1}(s)f(s^{-1})$, where
$\Delta$ is the modulus of the group.
 Let $\Sigma$ be the set of all (equivalence
classes of)  continuous unitary representations $\pi$ of $G$ in
Hilbert spaces $H_{\pi}$. For $f\in L^1(G)$, $\pi\in\Sigma$, we
put $\pi(f)=\int_G f(g)\pi(g)dg\in B(H_{\pi})$ with the integral
converging in the strong operator topology, and then
$$||f||_\Sigma=\sup_{\pi\in\Sigma}||\pi(f)||,$$
where $||\cdot||$ is the operator norm in $B(H_{\pi})$.

The {\it enveloping $C^*$-algebra} $C^*(G)$ of $G$ is the
completion  of $L^1(G)$ with respect to $||\cdot||_{\Sigma}$. Let
$\lambda:G\to B(L_2(G))$ be the left regular representation given
by $\lambda(s)f(g)=f(s^{-1}g)$. We denote by  $C^*_r(G)$ the {\it
reduced $C^*$-algebra} of $G$, that is the $C^*$-algebra generated
by operators $\lambda(f)\in B(L_2(G))$, $f\in L^1(G)$, and by
$VN(G)$ the von Neumann algebra of $G$, that is
$$VN(G)=\overline{\text{span} }^{WOT}\{\lambda(g):g\in G\}=
\overline{C^*_r(G)}^{WOT}\subset B(L_2(G)).$$

 The {\it Fourier-Stieltjes algebra}, $B(G)$, is the set of all coefficients
$s\mapsto \pi_{\xi,\eta}(s)=(\pi(s)\xi,\eta)$, where
$\pi\in\Sigma$, $\xi$, $\eta\in H_{\pi}$, of unitary
representations of $G$, as defined by Eymard, \cite{E}. $B(G)$ is
a Banach algebra with respect to the the norm
$$||u||=\text{inf}\{||\xi||\ ||\eta||:u=\pi_{\xi,\eta}\}.$$
Note that $B(G)\simeq C^*(G)^*$.

The {\it Fourier algebra}, $A(G)$,  is the family of functions
$s\mapsto (\lambda(s)\xi,\eta)=\bar\eta*\check\xi$, $\xi$, $\eta\in
L_2(G)$, $\check\xi(s)=\xi(s^{-1})$. $A(G)$ is identified with the
predual $VN(G)_*$ via
$\langle(\lambda(s)\xi,\eta),T\rangle=(T\xi,\eta)$, and thus is a
normed algebra with the norm denoted by $||\cdot||_A$. It is known
that for $u\in A(G)$ there exist even $\xi$, $\eta\in L_2(G)$ such
that $u=\bar\eta*\check\xi$ and $||u||_A=||\xi||_2\cdot ||\eta||_2$.
Furthermore, $A(G)$ is a closed ideal in $B(G)$.

Since $A(G)$ is the predual of the von Neumann algebra $VN(G)$,
$A(G)$ possesses a structure of operator space and one can define
a notion of completely bounded multipliers, $M_{cb}A(G)$, for
$A(G)$. For the theory of operator space and completely bounded
maps we refer the reader to \cite{effros-ruan,BS,nico}. A
complex-valued function $u:G\to {\mathbb C}$ is a completely
bounded multiplier if it is a multiplier, i.e. $uA(G)\subset
A(G)$, and is completely bounded as a linear map on $A(G)$. We
have $A(G)\subset B(G)\subset M_{cb}A(G)$ and, for $u\in A(G)$,
$\vert\vert u\vert\vert_{cb}\leq \vert\vert u\vert\vert_A$, where
$\vert\vert\cdot\vert\vert_{cb}$ is the completely bounded norm
(see \cite[Corollary 2.3.3]{nico}).
 We will use here a  characterization of
$M_{cb}A(G)$ obtained by N.Spronk in \cite{nico} as formulated  in
Theorem \ref{nico}.

\section{Spectral and operator synthesis}
Let $A$ be a semisimple, regular, commutative Banach algebra with
$X_A$ as spectrum; for any $a\in A$ we shall denote then by $\hat
a\in C_0(X_A)$ its Gelfand transform. Let also $E\subset X_A$ be a
closed subset. We then denote by
\begin{eqnarray*}
&I_A(E)= \{a \in A\mid \hat{a}^{-1}(0) \text{ contains } E\},\\
&J_A^0(E)=\{a\in A\mid \hat{a}^{-1}(0) \text{ contains a
neighborhood of }E\}
 \text{ and }
J_A(E)=\overline{J_A^0(E)}.
\end{eqnarray*}
It is known that $I_A(E)$ and $J_A(E)$ are  the largest and the
smallest closed ideals with $E$ as hull, i.e., if $I$ is a closed
ideal such that $\{x\in X_A: f(x)=0 \text{ for all } f\in I\}=E$
then $$ J_A(E)\subset I\subset I_A(E).$$ Let $I_A^c(E)$ denote the
set of all compactly supported functions $f\in I_A(E)$. We say
that $E$ is a {\it set of spectral synthesis} ({\it local spectral
synthesis}) for $A$ if $J_A(E)=I_A(E)$ ($I_A^c(E)\subset J_A(E)$).

Let $A^*$ be the dual of $A$. For $a\in A$ we set $\text{supp}(a)=
\overline{\{x\in X_A: \hat a(x)\ne 0\}}$ and
$\text{null}(a)=\{x\in X_A:\hat a(x)=0\}$. For $\tau\in A^*$ and
$a\in A$ define $a\tau$ in $A^*$ by $a\tau(b)=\tau(ab)$ and define
the support of $\tau$ by
$$\text{supp}(\tau)=\{x\in X_A: a\tau\ne 0\text{ whenever } a(x)\ne0\}.$$
Then, for
a closed set $E\subset X_A$
$$J_A(E)^{\perp}=\{\tau\in A^*:\text{supp}(\tau)\subset E\}$$
 and $E$ is spectral for $A$ if and only if
$\tau(a)=0$ for any $a\in A$ and $\tau\in A^*$ such that
$\text{supp}(\tau)\subset E\subset\text{null}(a)$.

The algebra $A(G)$ is a  semi-simple abelian regular Banach algebras
with spectrum $G$. In what follows we write $I_A(E)$ for
$I_{A(G)}(E)$, $J_A(E)$ for $J_{A(G)}(E)$,
 and
$\supp_{VN}T$ for $\supp (T)$ if $T\in VN(G)=(A(G))^*$.

Now we recall some definitions and important facts  on operator
synthesis following \cite{arv, sht}.  To make use of the results
from \cite{arv, sht} we will assume in the rest of the paper that
$G$ is a second countable locally compact group or s.c.l.c group
for short, and therefore is metrizable by
\cite[8.3]{hewitt-rossI}.

A subset $E\subset G\times G$ is called {\it marginally null}
(with respect to $m\times m$) if $E\subset (M\times G)\cup(G\times
N)$ and $m(M)=m(N)=0$. Two subsets $E_1$, $E_2$ are marginally
equivalent ($E_1\sim ^M E_2$ or simply $E_1\cong E_2$) if their
symmetric difference is marginally null. Furthermore,
$E_1\subset^M E_2$ means that $E_1\setminus E_2$ is marginally
null, a property holds marginally almost everywhere if it holds
everywhere apart of a marginally null set, and so on.

Following \cite{eks} we define a pseudo-topology on $G$. We call a
subset $E$ a pseudo-open if it is marginally equivalent to a
countable union of measurable rectangles $A\times B$. The
complements of pseudo-open sets are pseudo-closed sets.

Set  $T(G)=L^2(G)\hat\otimes L^2(G)$, where $\hat\otimes$ denotes
the projective tensor product. Note that in \cite{sht} the
notation $\Gamma(G,G)$ is used instead of $T(G)$. Every $\Psi\in
T(G)$ can be identified with a function $\Psi:G\times G\to{\mathbb
C}$ which admits a representation
\begin{equation}\label{eq1}
\Psi(x,y)=\sum_{n=1}^{\infty}f_n(x)g_n(y)
\end{equation}
where $f_n\in L^2(G)$, $g_n\in L^2(G)$   and
$\sum_{n=1}^{\infty}||f_n||_2\cdot||g_n||_2<\infty$.
 Such a representation
defines a function marginally almost everywhere (m.a.e.), so two
functions in  $T(G)$ which coincides m.a.e. are identified. The
$L^2(G)\hat\otimes L^2(G)$-norm of $\Psi$ is
$$||\Psi||_{T(G)}=\inf\{\sum_{n=1}^{\infty}||f_n||_2\cdot
||g_n||_2: \Psi=\sum_{n=1}^{\infty}f_n\otimes g_n\}.$$

 Note that a simple renormalization shows that each
$\Psi\in T(G)$ admits a representation (\ref{eq1}) such that
 $\sum_{n=1}^{\infty}||f_n||_2^2\cdot
\sum_{n=1}^{\infty}||g_n||_2^2<\infty$ and the norm
$||\Psi||_{T(G)}$ can be taken as the square root of the infinum
of $\sum_{n=1}^{\infty}||f_n||_2^2\cdot
\sum_{n=1}^{\infty}||g_n||_2^2$ over all such representations. By
\cite{eks} any $\Psi\in T(G)$ is pseudo-continuous. Thus if $\Psi$
vanishes (m.a.e.) on $K\subset G$ it vanishes on a pseudo-closed
set. For ${\mathcal F}\subset T(G)$, the null set, $\text{null
}{\mathcal F}$, is defined to be the largest pseudo-closed set
such that each function $F\in{\mathcal F}$ vanishes on it. For a
pseudo-closed set $E\subset G\times G$, let
\begin{eqnarray*}
&\Phi(E)=\{w\in T(G):w=0 \text{ m.a.e. on } E\},\\
&\Phi_0(E)=\overline{\{w\in T(G): w=0 \text{ on a
pseudo-neighbourhood of } E\}}.
\end{eqnarray*}
The spaces $\Phi(E)$, $\Phi_0(E)$ are $L^{\infty}(G)\times
L^{\infty}(G)$-bimodules: if $f, g\in L^{\infty}(G)$ and $w\in
\Phi(E)$ ($w\in \Phi_0(E)$) then $f(x)g(y)w(x,y)\in\Phi(E)$
($f(x)g(y)w(x,y)\in\Phi_0(E)$ respectively). Moreover, $\Phi(E)$
and $\Phi_0(E)$ are the largest and the smallest
$L^{\infty}(G)\times L^{\infty}(G)$-invariant subspaces of $T(G)$
whose null set is $E$.

A subset $E\subset G\times G$ is called a {\it set of (operator)
synthesis} or {\it synthetic} with respect to $m$ if
$\Phi(E)=\Phi_0(E)$.

It is know that $B(L_2(G))\simeq T(G)^*$ (see \cite{arv}). The
duality is given by
 $$\langle T,\Psi\rangle=\sum_{n=1}^{\infty}(Tf_n,\bar{g}_n),$$
for $T\in B(L_2(G))$ and $\Psi=\sum_{n=1}^{\infty}f_n\otimes g_n\in
T(G)$.

Let $P_U$ denote the multiplication operators by the characteristic
functions of a subset $U\subset G$. We say that $T\in B(L_2(G))$ is
{\it supported} in $E\subset G\times G$ (or $E$ supports $T$) if
$P_VTP_U=0$ for each pair of Borel sets $U\subset G$, $V\subset G$
such that $(U\times V)\cap E\cong \emptyset$.
 Then there exists the smallest (up to a marginally null set)
pseudo-closed set, $\supp T$, which supports $T$. We write
$\supp(T)\subset E$ if $T$ is supported by $E$. In the seminal paper
\cite{arv} Arveson defined a support in a similar way but using
closed sets instead of pseudo-closed. This closed support,
$\text{supp}_A T$ can be strictly larger than $\supp(T)$. Then $E$
is a set of operator synthesis if $\langle T,w\rangle=0$ for any
$T\in B(L_2(G))$ and $w\in T(G)$ with $\supp(T)\subset
E\subset\text{null }w$ (the inclusions up to a marginally null set).

We consider also the space $V^{\infty}(G)$ of all (marginal
equivalence classes of) functions $\Psi(x,y)$ that can be written in
the form (\ref{eq1}) with $f_n\in L^{\infty}(G)$, $g_n\in
L^{\infty}(G)$ and
$$\sum_{n=1}^{\infty} |f_n(x)|^2\leq C,\quad x\in G,\quad
\sum_{n=1}^{\infty}|g_n(y)|^2\leq C,\quad y\in G,$$ and for such
$\Psi$  we have
$$\vert\vert\Psi\vert\vert_{V^{\infty}}=\text{inf}\{\vert\vert
\sum_{n=1}^{\infty} |f_n|^2\vert\vert_{\infty}^{1/2}
\vert\vert\sum_{n=1}^{\infty}|g_n|^2\vert\vert^{1/2}_{\infty}:
\Psi=\sum_{n=1}^{\infty}f_n\otimes
g_n\}.$$ In tensor notations
$V^{\infty}(G)=L^{\infty}(G)\hat\otimes^{w^*h} L^{\infty}(G)$, the
weak$^*$-Haagerup tensor product (\cite{BS}). Let
$$V^{\infty}_{inv}(G)=\{w\in V^{\infty}(G): w(sr,tr)=w(s,t)\text{
for all $r$ in $G$ and m.a.e. $(s,t)\in G\times G$}\}.$$

In \cite{nico} Spronk found a connection between
$V_{inv}^{\infty}(G)$ and the algebra $M_{cb}A(G)$ of completely
bounded multipliers of $A(G)$. For a function $u: G\to {\mathbb
C}$ and  $t$, $s\in G$ define $$(Nu)(t,s)=u(ts^{-1}).$$
\begin{theorem}\label{nico}\cite{nico}
The map $u\mapsto Nu$ is a complete isometry from $M_{cb}A(G)$ onto
$V_{inv}^{\infty}(G)$.
\end{theorem}

\section{Spectral synthesis and operator synthesis}\label{sec3}

In this section we will prove our main result establishing a
connection between operator synthesis and spectral synthesis for
$A(G)$, where $G$ is a second countable locally compact group.
 The proofs are inspirited by the proof of
\cite[Theorem~4.6]{st}.

The Banach space  $B(L_2(G))$ is a left $V^{\infty}(G)$-module
with the action defined for $w=\sum_{i=1}^{\infty}\varphi_i\otimes
\psi_i\in V^{\infty}(G)$ and $T\in B(L_2(G))$ by
$$w\cdot T=\sum_{i=1}^{\infty}M_{\psi_i}TM_{\varphi_i},$$
where the partial sums converge strongly. The operator $T\mapsto
\sum_{i=1}^{\infty}M_{\psi_i}TM_{\varphi_i}$ we will also denote by
$\Delta_w$.

 For a closed subset $E\subset G$ we set
$$E^*=\{(s,t)\in G\times G\mid st^{-1}\in E\}.$$
\begin{lemma}\label{supp}
Let $S\in VN(G)$. Then
$$\supp (S)\subset \{\supp_{VN}S\}^*.$$
\end{lemma}
\begin{proof}
Let $U$, $V$ be closed  subsets of $G$ such that $(U\times
V)\cap\{\supp_{VN}(S)\}^*=\emptyset$. Then there exists an open
neighborhood, $W$, of $\{\supp_{VN}S\}^*$ such that $(U\times V)\cap
W=\emptyset$. Take $f$ and $g$ in $L_2(G)$ such that
$\supp(f)\subset U$, $\supp(g)\subset V$. For
$u=\langle\lambda(\cdot)f,g\rangle$, we have $(Sf,g)=\langle
S,u\rangle$. Moreover, $\supp(u)\subset \overline{UV^{-1}}$ and
$\supp(u)\cap\supp_{VN}S\subset
\overline{UV^{-1}}\cap\supp_{VN}S=\emptyset.$ Thus $0=\langle
S,u\rangle=(Sf,g)$. As $f$ and $g$ are chosen arbitrarily,
$P_VSP_U=0$. By the regularity of $m$, the last holds for any Borel
sets $U$, $V$ giving the statement.
\end{proof}

\begin{remark}\label{rem}\rm
Let $H$ be a closed subgroup of $G$. Then
\begin{eqnarray*}H^*&=&\{(s,t)\in G\times G: st^{-1}\in
H\}=\{(s,t)\in G\times G: Hs=Ht\}
\\&=&\{(s,t)\in G\times G:
f(t)=f(s)\}, \end{eqnarray*}

where $f:G\to H\backslash G$ is a continuous mapping defined by
$f(t)=Ht$.

Assume $\text{supp}_{VN}(S)\subset H$. By Lemma~\ref{supp},
$\text{supp}(S)\subset H^*$.  As for any Borel set $\Delta\subset
H\backslash G$ and $\alpha=f^{-1}(\Delta)$, we have
$(\alpha^c\times\alpha)\cap H^*=\emptyset$, it gives
$P_{\alpha^c}SP_{\alpha}=0$. Since this is true for any $\Delta$
we have also  $P_{\alpha}SP_{\alpha^c}=0$ implying that
$P_{\alpha} S=SP_{\alpha}$ and hence $S$ belongs to the commutant
${\mathcal B}'$ of the von Neumann algebra ${\mathcal B}$
generated by multiplication operators by functions $\varphi\in
L^{\infty}(G)$ which are constant on the right cosets.
\end{remark}
\begin{theorem}\label{synthesis}
Let $G$ be a s.c.l.c. group and $E\subset G$ be a closed subset.
If $E^*$ is synthetic with respect to $m$ then $E$ is a set of
local synthesis for $A(G)$.
\end{theorem}
\begin{proof}
Assume that $E^*$ is synthetic with respect to $m$. Let $u\in
I_A^c(E)$ and $S\in VN(G)$,
 $\supp_{VN}(S)\subset E$.
By Theorem~\ref{nico}, $Nu\in V^{\infty}(G)$. Moreover,
\begin{equation}\label{delta}
uT=\Delta_{Nu}T \text{ for any }T\in VN(G).
\end{equation}
In fact, if $Nu(s,t)=\sum_i\varphi_i(t)\psi_i(s)$, then
\begin{eqnarray*}
&&\Delta_{Nu}\lambda(s)f(t)=\sum_iM_{\varphi_i}\lambda(s)M_{\psi_i}f(t)=\\
&&=\sum_i\varphi_i(t)\psi_i(s^{-1}t)f(s^{-1}t)=Nu(t,s^{-1}t)f(s^{-1}t)=u(s)
\lambda(s)f(t)
\end{eqnarray*}
for any $f\in L^2(G)$ and $s\in G$.
The operator $\Delta_{Nu}$ is weakly-continuous. In fact, if
$S_k\to 0$ weakly, $||S_k||\leq C$ for some constant $C$ and
\begin{eqnarray*}
&&|(\Delta_{Nu}(S_k)f,g)|\leq\sum_{i=1}^n|(S_k\psi_if,\bar\varphi_ig)|+
\sum_{i=n+1}^{\infty}
|(S_k\psi_if,\bar\varphi_ig)|\\&&\leq\sum_{i=1}^n|(S_k\psi_if,\bar\varphi_ig)|+||S_k||
\left(\sum_{i=n+1}^{\infty}||\psi_if||^2\right)^{1/2}\left(\sum_{i=n+1}^{\infty}|
|\varphi_ig||^2\right)^{1/2}\\&&=
 \sum_{i=1}^n|(S_k\psi_if,\bar\varphi_ig)|+C
\left(\int_G\sum_{i=n+1}^{\infty}|\psi_i(t)f(t)|^2dt\right)^{1/2}
\left(\int_G\sum_{i=n+1}^{\infty}|\varphi_i(t)g(t)|^2dt\right)^{1/2}.
\end{eqnarray*}
 For given $\varepsilon>0$, by Lebesgue's theorem,  there
exists $n$ such that the second summand is less than $\varepsilon$
and then, as $S_k\to 0$ weakly, there exists $K$ such that the
first summand is less than $\varepsilon$ for any $k\geq K$.
Therefore (\ref{delta}) holds for any $T\in VN(G)$.

Clearly, since $u\in I_A(E)$,  $Nu$ vanishes on $E^*$. By
Lemma~\ref{supp}, we also have that $\supp(S)\subset
\{\text{supp}_{VN}S\}^*\subset E^*$. Therefore, for each $w\in
T(G)$,
$$\langle \Delta_{Nu}S, w\rangle=\langle S,(Nu)w\rangle=0$$
so that  $uS=\Delta_{Nu}S=0$.

 From the regularity of $A(G)$ it
follows that there exists a compactly supported function $v\in A(G)$
such that $v=1$ on the support of $u$. Thus
$$\langle S,u\rangle=\langle S,vu\rangle=\langle uS,v\rangle=0.$$
\end{proof}
It is easy to see that the condition for $E$ to be a set of local
synthesis for $A(G)$ is equivalent to the condition $uS=0$ for any
$u\in A(G)$ and $S\in VN(G)$ such that $\supp_{VN}(S)\subset
E\subset \text{null } u$. If $G$ is amenable then we have the
implication
\begin{equation}\label{H}
uT=0\Rightarrow \langle T,u\rangle=0,
\end{equation}
guaranteed by the existence of a bounded approximate identity
 and any set of local synthesis is a set of synthesis.
Certainly the assumption of boundedness of the identity is
superfluous, and the statement holds even for $G$ such that $A(G)$
has an unbounded approximate identity. So that we have
\begin{cor}\label{ai}
Let $G$ be  a s.c.l.c.\ group such that $A(G)$ has an (unbounded)
approximate identity  and let $E$ be a closed subset of $G$. If
$E^*$ is a set of operator synthesis with respect to $m$ then $E$
is a set of spectral synthesis.
\end{cor}
It is not known whether an approximate identity  exists for any
locally compact group $G$. Unbounded approximate identities which
are completely bounded as multipliers of the Fourier algebra $A(G)$
were studied in \cite{can_haag,can,cowhaag}. Those exist for a
number of groups like the general Lorentz group $SO_0(n,1)$, its
closed subgroups, in particular, the free group ${\mathbb F}_n$ on
$n$ generators, extensions of $SO_0(n,1)$ by a finite group, in
particular, $SL(2,{\mathbb R})$, $SL(2,{\mathbb C})$, $SL(2,{\mathbb
H})$, and weakly-amenable groups.

The property (\ref{H}) of an operator $T\in VN(G)$ was discussed
in \cite{E} and called there by the property $(H)$. In particular,
it was shown that any $T$ supported in a compact set $E$ possesses
this property. Thus any compact set of local synthesis is a set of
spectral synthesis.
\begin{cor}\label{compact}
Let $E$ be a compact subset of $G$. If $E^*$ is a set of synthesis
with respect to $m$ then $E$ is a set of spectral synthesis.
\end{cor}
It is  not known whether there exists $T\in VN(G)$ which does
not satisfy (\ref{H}).

The next statement was proved by Eymard for arbitrary locally
compact groups $G$, \cite{E}, using more complicated arguments of
the theory of distributions on  $G$.
\begin{cor}
 Let $s$ be an element of the s.c.l.c. group $G$. Then
$\{s\}$ is a set of spectral synthesis.
\end{cor}
\begin{proof} It is enough to prove the statement  for
$E=\{e\}$, where $e$ is the identity element in $G$. In this case
$E^*=\{(s,s):s\in G\}$. That $E^*$ is a set of operator synthesis
follows e.g. from \cite[Theorem~4.8]{sht}, but it can be easily seen
also using the following simple arguments.

Let $T\in B(L_2(G))$ and $\supp (T)\subset E^*$. It follows from
Remark~\ref{rem} that $T$ is the multiplication operator by some
function $a\in L^{\infty}(G,m)$. Let
 now $F=\sum_{i=1}^{\infty}f_i\otimes g_i\in \Phi(E^*)$. Then
 $$\langle T,F\rangle=\sum_{i=1}^{\infty}(Tg_i,\bar f_i)=
 \sum_{i=1}^{\infty}(ag_i,\bar f_i)=\int_Ga(r)F(r,r)dr=0,$$

and therefore $E^*$ is a set of synthesis with respect to $m$.
 By Corollary~\ref{compact}, $E$ is a set of spectral synthesis.
\end{proof}
We also have another proof in the case of s.c.l.c.\ groups of the
following statement proved by Herz, \cite[Theorem~2]{herz}.
\begin{cor}
Any closed subgroup $H$ of $G$ is a set of local spectral synthesis.
\end{cor}
\begin{proof}
We have $H^*=\{(s,t)\in G\times G: st^{-1}\in H\}=\{(s,t)\in
G\times G: Hs=Ht\}=\{(s,t)\in G\times G: f(t)=f(s)\}$, where
$f:G\to H\backslash G$ is a continuous function defined by
$f(t)=Ht$. By \cite[Theorem~4.8]{sht}, $H^*$ is a set of operator
synthesis and the statement now follows from
Theorem~\ref{synthesis}.
\end{proof}

\begin{remark}\rm
By Corollaries~\ref{ai}, \ref{compact}, any compact subgroup $H$ is
spectral for $A(G)$ and any closed subgroup is spectral for $A(G)$
if $A(G)$ has an approximate identity.

That any closed subgroup of a locally compact group $G$ is a set
of spectral synthesis was shown by Takesaki and Tatsuuma,
\cite[Theorem~3]{taketat}, using the Mackey imprimitivity theorem
and the result mentioned in Remark~\ref{rem} proved by a different
method.
\end{remark}
\begin{example}\rm
Let $R$ be an ordered s.c.l.c.\ group and $\Delta:G\to R$ be a
continuous homomorphism. Take a finite intersection of intervals
$\cap_{k=1}^{m}[\alpha_k, \beta_k]$ in $R$ and set
$E=\Delta^{-1}(\cap_{k=1}^{m}[\alpha_k, \beta_k])$. Then
\begin{eqnarray*}
  E^* &=&\{(s,t):\alpha_k\leq\Delta(st^{-1})\leq\beta_k,k=1,\ldots,m\}\\
  &=&\{(s,t):\alpha_k\Delta(t)\leq\Delta(s)
\leq\beta_k\Delta(t), k=1,\ldots,m\}\\
  &=& \{(s,t):f_i^k(t)\leq g_i^k(s), i=1,2, k=1,\ldots,m\},
  \end{eqnarray*}

where $f_1^k(t)=\alpha_k\Delta(t)$, $f_2^k(t)=\Delta(t)^{-1}$,
$g_1^k(s)=\Delta(s)$, $g_2^k(s)=\Delta(s)^{-1}\beta_k$. By
\cite[Theorem~4.8]{sht}, $E^*$ is a set of operator synthesis with
respect to the Haar measure on $G$. Therefore $E$ is a set of
spectral synthesis by Theorem~\ref{synthesis}.
\end{example}
Our next aim is to prove a   converse to the statement of
Theorem~\ref{synthesis}.

The Banach space $T(G)$ is an $L^1(G)$-module with the action
defined by
$$f\odot w(s,t)=\int_G f(r)\Delta^{1/2}(r)w(sr,tr)dr, \quad f\in L^1(G), w\in T(G),$$
where $\Delta$ is the modular function of $G$. We have $||f\odot
w||_{T(G)}\leq||f||_1\cdot||w||_{T(G)}$. Moreover, $e_{\alpha}\odot
w\to w$ for any bounded approximate identity $\{e_{\alpha}\}$ in
$L^1(G)$ (see \cite[32.22, 32.33]{hewitt-rossII} and
\cite[p.365]{st}).

Let us define another action by compactly supported
$L^1(G)$-functions $f$:
$$f\cdot w(s,t)=\int_G f(r)w(sr,tr)dr, \quad  w\in T(G).$$

Then $f\odot w=f\Delta^{1/2}\cdot w$.\\

We observe that by the estimate (\ref{l1estim}) below, the
integral $\int_G f(r)w(sr,tr)dr,  w\in T(G), s,t\in G,$ converges
also if $f\in L^{\infty}(G)$ and so defines a mapping $(s,t)\to
f\cdot w(s,t):=\int_G f(r)w(sr,tr)dr.$

\begin{theorem}\label{ops}
Let $G$ be a s.c.l.c.\ group. If a closed subset $E\subset G$ is
a set of local spectral synthesis for $A(G)$ then  $E^*$ is synthetic with respect
to Haar
measure.
\end{theorem}

\begin{proof}

It is sufficient to show that $w\cdot T=0$ for $T\in B(L_2(G))$ and
$w\in V^{\infty}(G)$ such that $\supp (T)\subset
E^*\subset\text{null }w$ (\cite[Proposition~5.3]{sht2}).

As $G$ is second countable the group $G$ is $\sigma$-compact and
therefore there exist compact sets $K_n$ such that $K_n\subset
K_{n+1}$ and $\cup_{n=1}^{\infty}K_n=G$. Then clearly,
$M_{\chi_{K_n}}TM_{\chi_{K_n}}\to T$ strongly. Therefore we can
restrict ourselves to a compactly supported operator $T$, i.e.
$\supp (T)\subset M\times M$ for a compact set $M\subset G$, and
compactly supported $w\in V^{\infty}(G)$. Note that in this case
$w\in T(G)$.

Let $\hat G$ denote the set of (equivalence classes of)
irreducible continuous unitary representations of $G$. For
$\pi\in\hat G$ let $H_{\pi}$ denote the representation space of
$\pi$. Fix a basis $\{e_j\}_{j}$ in $H_{\pi}$ and denote by
$u_{kj}^{\pi}$ the matrix coefficients of $\pi$, i.e.
$u_{kj}^{\pi}(s)=(\pi(s)e_k,e_j)$.

For (a compactly supported)
$w=\sum_{i=1}^{\infty}\varphi_i\otimes\psi_i\in T(G)$ and
$\pi\in\hat G$ consider now the following operator-valued function
$$w^{\pi}(s,t)=\int_G w(sr,tr)\pi(r)dr.$$
The integral is well-defined as a Bochner integral. In fact, for
each $s$, $t\in G$, applying Cauchy-Schwartz's inequality, we obtain
\begin{eqnarray}\label{l1estim}
\int_G|w(sr,tr)|dr&\leq&\int_G\sum_{i=1}^{\infty}|\varphi_i(sr)\psi_i(tr)|dr\\
\nonumber&\leq&\int_G\left(\sum_{i=1}^{\infty}|\varphi_i(sr)|^2\right)^{1/2}
\left(\sum_{i=1}^{\infty}|\psi_i(tr)|^2\right)^{1/2}dr\\
\nonumber&\leq&\left(\sum_{i=1}^{\infty}\int_G|\varphi_i(sr)|^2dr\right)^{1/2}
\left(\sum_{i=1}^{\infty}\int_G|\psi_i(tr)|^2dr\right)^{1/2}\\
\nonumber&=&
\left(\sum_{i=1}^{\infty}||\varphi_i||_2^2\right)^{1/2}
\left(\sum_{i=1}^{\infty}||\psi_i||_2^2\right)^{1/2}<\infty.
\end{eqnarray}
Set $\tilde w^{\pi}(s,t)=\pi(s)w^{\pi}(s,t)$ and
\begin{equation}\label{wkj}w_{kj}^{\pi}(s,t)=(w^{\pi}(s,t)e_k,e_j)=u_{kj}^{\pi}\cdot
  w(s,t),
\quad
\tilde w_{kj}^{\pi}(s,t)=(\tilde
w^{\pi}(s,t)e_k,e_j).\end{equation}

If $w\in\Phi(E^*)$ then $w(sr,tr)$ vanishes m.a.e.\  on $E^*$ for
all $r$, and therefore  $w^{\pi}(s,t)$, $\tilde w^{\pi}(s,t)$,
$w^{\pi}_{kj}(s,t)$ and $\tilde w^{\pi}_{kj}(s,t)$ vanish on
$E^*$.

We have the following expression for $w_{kj}^\pi$ and $\tilde
w_{kj}^\pi$:
\begin{eqnarray*}
\tilde w_{kj}^{\pi}(s,t)=\int_G
w(sr,tr)(\pi(sr)e_k,e_j)dr&=&\int_G w(r,ts^{-1}r)
(\pi(r)e_k,e_j)dr\\&=&\int_G w(r,ts^{-1}r) u_{kj}^{\pi}(r)dr.
\end{eqnarray*}
 In
particular, if $w=f\otimes g\in T(G)$, then
\begin{equation}\label{nfukl} \tilde
w_{kj}^{\pi}(s,t)=\int_G f(r)u_{kj}^\pi(r)g(ts\inv
r)dr=N(fu_{kj}^\pi*\check g)(s,t).\end{equation}

Furthermore
\begin{eqnarray}\label{wkjpi}
w_{kj}^{\pi}(s,t)&=&\int_G w(sr,tr)(\pi(r)e_k,e_j)dr=\int_G
w(r,ts^{-1}r)
(\pi(s^{-1}r)e_k,e_j)dr\nonumber\\
&=&\int_G w(r,ts^{-1}r) (\pi(r)e_k,\pi(s)e_j)dr\nonumber\\&=&\sum_{l}\int_G
w(r,ts^{-1}r)
(\pi(r)e_k,e_l)(e_l,\pi(s)e_j)dr\nonumber\\
&=&\sum_{l} u_{lj}^{\pi}(s)\int_G w(r,ts^{-1}r)
u_{kl}^{\pi}(r)dr=\sum_lu_{lj}^{\pi}(s)\tilde w_{kl}^{\pi}(s,t).
\end{eqnarray}

We state first that  $\tilde w^{\pi}_{kj}(s,t)\in V^{\infty}(G)$.
Indeed, if $w=f\otimes g\in T(G)$, then by (\ref{nfukl}) $\tilde
w_{kj}^\pi=N(fu_{kj}^\pi\ast\check g)$ and
 and therefore, since
$fu_{kj}^{\pi}*\check g\in A(G)$, $\tilde w_{kj}^\pi\in
V^{\infty}(G)$, by Theorem~\ref{nico}. Moreover,
$$||\tilde w_{kj}^{\pi}(s,t)||_{V^{\infty}}\leq ||fu_{kj}||_2||g||_2\leq
||f||_2||g||_2,$$

 so that the linear operator $w\mapsto \tilde
w_{kj}^{\pi}\in V^{\infty}(G)$  defined on elementary tensors
extends to a bounded operator $T(G)\to V^{\infty}(G)$. Thus
$\tilde w_{kj}^{\pi}\in V^{\infty}(G)$ for any $w\in T(G)$. \\

Next we show that $w_{kj}^\pi\in V^\infty (G)$. For $T\in
B(L^2(G))$ and $w=f\otimes g\in T(G)$ such that $\no{w}=\no f _2
\no g_2$, define $w_{kj}^{\pi}\cdot T$ by
\begin{equation}\label{wkjui}\langle w_{kj}^{\pi}\cdot T,\Psi\rangle=
\sum_{l}\langle\tilde w_{kl}^{\pi}\cdot
T,u_{lj}^{\pi}(s)\Psi(s,t)\rangle,\quad \Psi\in
T(G).\end{equation}

 This formula makes sense. In fact, if
$\Psi=\sum_{i=1}^{\infty} f_i\otimes g_i$, such that $\no \Psi^2
_{T(G)} =\sum_{i=1}^\infty \no {f_i}^2\sum_{i=1}^\infty \no
{g_i}^2$, we have

\begin{eqnarray*}\vert\sum_{l}\langle\tilde w_{kl}^{\pi}\cdot
T,u_{lj}^{\pi}\Psi\rangle\vert^2 &\leq &
\left(\sum_{l}\vert\vert\tilde w_{kl}^{\pi}\cdot
T\vert\vert^2\right)\left(\sum_{l}\vert\vert
u_{lj}^{\pi}\Psi\vert\vert_{T(G)}^2\right)\\
&\leq & \left(\sum_{l}\vert\vert\tilde
w_{kl}^{\pi}\vert\vert_{V^{\infty}}^2\vert\vert
T\vert\vert^2\right)\left(\sum_{l}\sum_{i=1}^{\infty}\vert\vert
u_{lj}^{\pi}f_i\vert\vert_{2}^2\cdot\sum_{i=1}^{\infty}\vert\vert
g_i\vert\vert_{2}^2\right)\\&\leq &\vert\vert
T\vert\vert^2\left(\sum_{l}
||fu_{kl}^\pi||_2^2||g||_2^2\right)\left(\sum_{i=1}^{\infty}\sum_{l}\vert\vert
u_{lj}^{\pi}f_i\vert\vert_{2}^2\right)\cdot\sum_{i=1}^{\infty}\vert\vert
g_i\vert\vert_{2}^2\\
&=&\vert\vert
T\vert\vert^2||f||_2^2||g||_2^2\sum_{i=1}^{\infty}||f_i||_2^2\cdot\sum_{i=1}^{\infty}|
|g_i||_2^2\\&\leq&
\vert\vert T\vert\vert^2||w||_{T(G)}^2||\Psi||_{T(G)}^2,
\end{eqnarray*}
the last equality follows from
\begin{eqnarray*}\sum_{l} ||fu_{kl}||_2^2&=&
\sum_{l}\int_G|f(t)u_{kl}^\pi(t)|^2dt=
\int_G|f(t)|^2\sum_{l}|(\pi(t)e_k,e_l)|^2dt\\
&=&\int_G|f(t)|^2||(\pi(t)e_k||_2^2dt= \int_G|f(t)|^2dt=||f||_2^2.
\end{eqnarray*} Thus the operator
$w_{kj}^{\pi}\cdot T$ is well-defined for $w=f\otimes g$ and since
$\vert\vert w_{kj}^{\pi}\cdot T\vert\vert \leq\vert\vert
T\vert\vert||w||_{T(G)}$ for elementary tensors $w$, the
definition $w_{kj}^{\pi}\cdot T$ makes sense for any $w\in T(G)$
and
\begin{equation}\label{inf}
\vert\sum_{l}\langle\tilde w_{kl}^{\pi}\cdot
T,u_{lj}^{\pi}\Psi\rangle\vert^2\leq\vert\vert
T\vert\vert^2||w||_{T(G)}^2||\Psi||^2_{T(G)}
\end{equation}
for any $\Psi\in T(G)$
Clearly, $T\mapsto w_{kj}^{\pi}\cdot T$ is thus a bounded
$L^{\infty}(G)$-bimodule map on $B(L^2(G))$ and hence by
\cite[2.1]{smith} is completely bounded. Then by \cite[4.2]{BS},
it is of the form $\omega\cdot T$ for $\omega\in V^{\infty}(G)$
and therefore $w_{kj}^{\pi}=\omega$ m.a.e.\\

For $\tilde w^{\pi}$ we
have
\begin{eqnarray*}
\tilde w^{\pi}(sr,tr)&=&\pi(sr)w^{\pi}(sr,tr)=\pi(sr)\int_G
w(srp,trp)
\pi(p)dp\\
&=&\pi(s)\int_G w(srp,trp) \pi(rp)dp=\pi(s)\int_G w(sp,tp)
\pi(p)dp=\tilde w^{\pi}(s,t),
\end{eqnarray*}
implying $\tilde w_{kj}^{\pi}\in V^{\infty}_{inv}(G)$. By
Theorem~\ref{nico},
 \begin{equation}\label{equal}\tilde
w_{kj}^{\pi}=Nu\end{equation}
 for some $u\in M_{cb}A(G)$.
Moreover, if $w$ vanishes on $E^*$ then $u$ vanishes on $E$.

 We claim that  $Nu\cdot
T=0$ for any operator $T$ and $u\in M_{cb}A(G)$ such that  $\supp
(T)\subset E^*$, $\text{null }u\supset E$. In fact, since $E$ is a
set of local spectral synthesis,  given $w\in A(G)$ with compact support, $uw$
 can be approximated by $u_{\alpha}\in J_A^0(E)$ and therefore $N(uw)$ is a
 $V^{\infty}(G)$-limit of
$Nu_{\alpha}$, vanishing on pseudo-neighborhoods of $E^*$. By
\cite[Theorem~4.3]{sht}, $\langle T,\Psi\rangle=0$ for any
$\Psi\in\Phi_0(E^*)$. Therefore $\langle Nu_{\alpha}\cdot
T,F\rangle=\langle T, (Nu_{\alpha})F\rangle=0$ for any $F\in T(G)$
implying $Nu_{\alpha}\cdot T=0$ and
 $N(uw)\cdot T=0$. As $\langle Nu\cdot
T,(Nw)F\rangle=0$ for any $w\in A(G)$ with compact support and $F\in
T(G)$ it is enough to see now that the subspace, ${\mathcal M}$,
generated by $(Nw)F$, where $F\in T(G)$ and $w\in A(G)$ with compact
support, is dense in $T(G)$. As $\text{null }{\mathcal
M}=\emptyset$, this is true by   an analogue of Wiener's Tauberian
Theorem \cite[Corollary~4.3]{sht}.

 We
obtain by (\ref{equal}) that for any $w\in V^{\infty}(G)$ which is
compactly supported and vanishes on $E^*$
$$\tilde w_{kj}^{\pi}\cdot T=0.$$

Therefore by (\ref{inf})
\begin{eqnarray*}
w_{kj}^{\pi}\cdot T =(\sum_{l}u_{jl}^{\pi}(s)\tilde
w^{\pi}_{kj}(s,t))\cdot T=0,
\end{eqnarray*}

Let $K$ be a compact set such that $\supp(w)\subset K\times K$ and
$\supp(T)\subset K\times K$. Then since $T=M_{\chi_K}TM_{\chi_K}$,
for compactly supported $h\in L^1(G)$ we have
\begin{eqnarray*}
\langle T,h\cdot w\rangle&=&\langle
M_{\chi_K}TM_{\chi_K}, h\cdot w\rangle=\sp{T}{\chi_K(t)\chi_K(s)\int_G w(sr,tr)h(r)dr}\\
&=&\langle T,\int_Gw(sr,tr)\chi_{K^{-1}K}(r)h(r)dr\rangle= \langle
T,h\chi_{K^{-1}K}\cdot w\rangle.
\end{eqnarray*}
Take  $f\in L^{\infty}(G)$, $g\in L^{\infty}(G)$, $\supp(f)\subset
K$, $\supp(g)\subset K$, and set $\omega=f\otimes g$. Then
\begin{eqnarray*} \langle \omega\cdot T,u_{kj}^\pi
\chi_{K^{-1}K}\cdot w\rangle= \langle
T,(u_{kj}^\pi\chi_{K^{-1}K}\cdot w)
\omega\rangle\\
=\langle T,w_{kj}^\pi\omega\rangle= \langle w_{kj}^{\pi}\cdot
T,\omega\rangle=0. \end{eqnarray*} Hence for finite linear
combinations $\sum_i c_i u_i$ of matrix coefficients
$$\langle\omega\cdot T,\sum_i c_iu_i\chi_{K^{-1}K}\cdot w\rangle=0.$$

 Take now an approximate identity $\{e_{\alpha}\}$ consisting
of non-negative  continuous functions with compact support in
$L^1(G)$. We can assume that $K^{-1}K$ contains the supports of
$e_{\alpha}$'s. Then, by \cite[13.6.5]{di}, $\DE^{1/2}e_{\alpha}$
 can be approximated in $L^1(G)$ by finite linear
combinations $\sum_{i=1}^d c_iu_i\chi_{K\inv K}$. This yields
$$0=\langle \omega\cdot T, (\DE^{1/2}e_{\alpha})\cdot w\rangle=
\langle \omega\cdot T,(e_{\alpha})\odot w\rangle$$
 and therefore $0=\langle \omega\cdot T,w\rangle=
 \langle w\cdot T,\omega\rangle$. Finally, we obtain
 $w\cdot T=0$.

\end{proof}
\begin{cor}
Let $G$ be a s.c.l.c.\  group.

a) Then a compact set $E\subset G$ is a set of spectral synthesis
for $A(G)$ if and only if $E^*$ is a set of operator synthesis with
respect to  Haar measure.

b) Assume that $A(G)$ has an approximate identity. Then a closed set
$E\subset G$ is a set of spectral synthesis for $A(G)$ if and only
if $E^*$ is a set of operator synthesis with respect to   Haar
measure.
\end{cor}
\begin{proof}
Follows from Corollary~\ref{ai}, Corollary~\ref{compact} and
Theorem~\ref{ops}.
\end{proof}

\section{Ditkin sets and operator Ditkin sets.}
Our goal in this section is to find a connection between Ditkin
sets and operator Ditkin sets. Let $G$ be a locally compact group
and let $m$ be the Haar measure on $G$. A closed subset $E\subset
G$ is said to be a (local) Ditkin set if for any $f\in I_A(E)$
($f\in I_A^c(E)$) there exists a sequence
$\{\tau_n\}_{n=1}^{\infty}$ such that $\tau_n\in J_A^0(E)$
($n=1,2,\ldots$) and $\tau_nf\to f$ as $n\to\infty$. It is called
a strongly Ditkin set if such a sequence $\{\tau_n\}$ can be
chosen uniformly for all functions $f\in I_A(E)$.

In a similar way operator Ditkin sets were defined in \cite{sht}. By
\cite{st} the elements of $V^{\infty}(G)$ are the multipliers of
$T(G)$, i.e. $w\omega\subset\omega$ if $w\in V^{\infty}(G)$ and
$\omega\in T(G)$. We call a pseudo-closed subset $E\subset G\times
G$ an $m$-Ditkin set if for any $w\in\Phi(E)$ there exists a
sequence $\tau_n\in V^{\infty}(G)$ such that $\tau_n$ vanishes on a
pseudo-neighbourhood of $E$ ($n=1,2,\ldots$) and
$$\vert\vert \tau_nw-w\vert\vert_{T(G)}\to 0,\text{ as
}n\to\infty.$$

$E\subset G\times G$ is said to be a strong $m$-Ditkin set if such
a sequence $\{\tau_n\}_{n=1}^{\infty}$ can be chosen uniformly for
all $w\in T(G)$.

If $G$ is a compact metrizable abelian group, it is known that
$E=\{0_G\}$ is a strong Ditkin set. Moreover, this set  satisfies
the following three conditions: (1) $E$ is a set of synthesis; (2)
there exist open sets $\Omega_n$ containing $E$ such that
$$\Omega_{n+1}\subset\Omega_n\quad n=1,2,\ldots\text{ and
}\cap_{n=1}^{\infty}\bar\Omega_n=E;$$ (3) there exists a sequence
$\{u_n\}$ with $1-u_n\in J_A^0(E)$, $n=1,2,\ldots$, satisfying the
following two conditions: \begin{eqnarray*} &u_n(x)=0\text{ for all
}x\notin\Omega_n,\\
&||u_n||\leq 1+\varepsilon_n,
\end{eqnarray*}
where $\{\varepsilon_n\}$ is a sequence decreasing to zero.

Note that  any closed subset $E$ of the spectrum of  a semisimple,
regular, commutative Banach algebra $A$ satisfying the conditions
(1), (2) and (3) is a strong Ditkin set for $A$ (\cite{drury}).

 Let $Z$ be a
standard Borel space and let $f:G\to Z$ and $g:G\to Z$ be Borel
functions. Consider
\begin{equation}\label{nontrian}E=\{(s,t):
f(s)=g(t)\}\subset G\times G.\end{equation}

 If $G$, $Z$ are
compact metrizable spaces, $f$, $g$ are continuous functions and
 $A=V(G)=C(G)\hat\otimes C(G)$, the Varopoulos algebra,
then $E$ is shown in \cite{drury} to satisfy the conditions (1), (2)
and (3) and therefore to be a strong Ditkin set for  $V(G)$. Knowing
that $E$ is a set of synthesis with respect to $m$
(\cite[Theorem~4.8]{sht}) we can use a similar argument to show the
following statement.

\begin{proposition}\label{nontr}
$E=\{(s,t): f(s)=g(t)\}\subset G\times G$ is a strong $m$-Ditkin
set.
\end{proposition}
\begin{proof}  First we embed $Z$ into the torus, ${\mathbb
T}_{\infty}$, of infinite dimension: by \cite[Theorem A.1]{ta}
there exists a Borel injective mapping $\psi:Z\to{\mathbb
T}_{\infty}$. Consider a mapping $\rho:G\times G\to {\mathbb
T}_{\infty}$ given by
$$\rho(s,t)=\psi(f(t))-\psi(g(s)),$$
where the subtraction is taken with respect to the group structure
on ${\mathbb T}_{\infty}$.
Then $\rho^{-1}({\mathbb T}_{\infty})=E$.

As $\{0_{\mathbb T_{\infty}}\}$ satisfies $(1)$, $(2)$, $(3)$ in
$A({\mathbb T}_{\infty})$, there exist open sets
$\Sigma_n\subset{\mathbb T}_{\infty}$ such that $0_{{\mathbb
    T}_{\infty}}\in \Sigma_n$,
$\Sigma_n\subset\Sigma_{n+1}$, $n=1,2\ldots$,
$\bigcap\bar\Sigma_n=\{0\}$ and a sequence of functions $w_n\in
A({\mathbb T}_{\infty})$ such that $1-w_n\in J_A^0(\{0_{\mathbb
T_\infty}\})$, $w_n(x)=0$ for $x\notin\Sigma_n$ and $\vert\vert
w_n\vert\vert_A\leq 1+\varepsilon_n$ ($\varepsilon_n\to 0$,
$\varepsilon_n>0$).

Define  $\Omega_n=\rho^{-1}(\Sigma_n)$. As $\psi\circ f$ and
$\psi\circ g$ are  Borel mappings for given $m>0$ there exist by
Lusin's theorem  closed subsets $A_m\subset G$ and $B_m\subset G$
such that $\rho$ is continuous on $A_m\times B_m$ and $|G\setminus
A_m|<1/m$, $|G\setminus B_m|<1/m$. We can choose those subsets
increasing in $m$. Therefore
$$\Omega_n\cap(A_m\times B_{m})=\{(x,y)\in A_{m}\times
B_{m}:\rho(x,y)\in\Sigma_n\}$$
 is open in $A_{m}\times B_{m}$ for all $n$ and
\begin{eqnarray}\label{inclusion}
E\cap (A_{m}\times
B_{m})&=&\cap_{n=1}^{\infty}\rho^{-1}(\Sigma_n)\cap(A_{m}\times
B_{m})\\
\nonumber &\subset&{\cap_{n=1}^{\infty}}
\overline{\Omega_n\cap(A_{m}\times B_{m})}\subset\\
&\subset &\cap_{n=1}^{\infty}\rho^{-1}(\overline{\Sigma_n})\cap
(A_{m}\times B_{m})=E\cap(A_{m}\times B_{m})\nonumber
\end{eqnarray}

 Set
$u_n=w_n\circ\rho$. First we show that $u_n\in V^{\infty}(G)$. In
fact, if $w_n=\sum_{\chi\in \hat{\mathbb
T}_{\infty}}a_{\chi,n}\chi$ with $\sum_{\chi\in\hat {\mathbb
T}_{\infty}}\vert a_{\chi,n}\vert=\no {w_n}_A$, then
$$w_n(\rho(s,t))=\sum_{\chi\in\hat {\mathbb
T}_{\infty}}a_{\chi,n}\chi(\psi(f(s)))\overline{\chi(\psi(g(t)))}$$

and our claim follows since  $$\sum_{\chi\in{\mathbb
T}_{\infty}}\vert a_{\chi,n}\vert\vert
\chi(\psi(f(s)))\vert^2=\sum_{\chi\in\hat {\mathbb
T}_{\infty}}\vert a_{\chi,n}\vert\vert
\chi(\psi(g(t)))\vert^2=\sum_{\chi\in\hat {\mathbb
T}_{\infty}}\vert a_{\chi,n}\vert=\vert\vert w_n\vert\vert_A.$$
 Moreover, by Theorem~\ref{nico}, $\vert\vert
u_n\vert\vert_{V^{\infty}}\leq\vert\vert w_n\vert\vert_A\leq
1+\va_n$, $u_n=0$ on $\Omega_n^c$ and $\tau_n=1-u_n$ vanishes on a
pseudo-neighbourhood  of $E$.

We next  show that for given $w\in\Phi(E)$, $||\tau_n w-w||\to 0$,
as $n\to 0$. Assume first that $\text{supp }(w)\subset K\times K$,
where $K$ is a compact set. Then $w=w_1^m+w_2^m+w_3^m+w_4^m$, where
$\displaystyle w_1^m=w\chi_{A_{m}\times
  B_{m}}$, $\displaystyle w_2^m=w\chi_{(G\setminus A_{m})\times B_{m}}$
   $\displaystyle w_3^m=w\chi_{A_{m}\times
  (G\setminus B_{m})}$, $\displaystyle w_4^m=w\chi_{(G\setminus A_{m})
  \times (G\setminus
  B_{m})}$. For given $\va>0$ there exists $M>0$ such that
  $||w_i^m||_{T(G)}<\va$ for each $m\geq M$ and $i=2,3,4$. Indeed, since the
  sequence $\{A_m\}$ is increasing in measure to $G$, by Lebesgue's
  theorem $\int_{G\setminus
  A_m}\sum_{i=1}^{\infty}|f_i(r)|^2 dr\to 0$ as $m\to\infty$
  and for
  $ w_2^m$ we have
  $$||w_2^m||^2_{T(G)}\leq \int_{G\setminus
  A_m}\sum_{i=1}^{\infty}|f_i(r)|^2 dr\int_{G}\sum_{i=1}^{\infty}|g_i(r)|^2
  dr\to 0.$$
Similarly, $||w_i^m||^2_{T(G)}\to 0$, $i=3,4$. Fix now $m>M$.
  By \cite[Theorem~4.8]{sht}  $E$  is a set of operator synthesis
with respect to Haar measure and by \cite[Lemma~6.1]{sht} so is
$E\cap(A_{m}\times B_{m})$. As $E\cap(A_m\times B_m)$ is closed,
by \cite[2.29]{arv} there exists $\psi\in T(G)$, $\text{supp
}\psi\subset A_{m}\times B_{m}$, vanishes on an open neighborhood
$\Omega\subset A_{m}\times B_{m}$ of $E\cap( A_{m}\times B_{m})$
and $||w_1^m-\psi||<\va$. By (\ref{inclusion}),
$$\cap_{n=1}^{\infty}
\overline{\Omega_n\cap(A_{m}\times B_{m})}\cap K\times
K\subset\Omega\cap K\times K$$

and so $\Omega_n\cap(A_{m}\times B_{m})\cap K\times
K\subset\Omega\cap K\times K$ for $n>N$, $N$ is large enough. Thus
$u_n\psi\chi_{K\times K}=0$ for $n>N$. We get
$$\vert\vert\tau_nw_1^m-w_1^m\vert\vert=\vert\vert
\tau_n(w_1^m-\psi)\chi_{K\times K}-u_n\psi\chi_{K\times
  K}-(w_1^m-\psi)\chi_{K\times K}\vert\vert\leq
(1+||\tau_n||_{V^{\infty}})\va$$ and
$$\vert\vert\tau_nw-w\vert\vert\leq
\vert\vert\tau_nw_1^m-w_1^m\vert\vert+
\sum_{i=2}^4\vert\vert(\tau_n-1)w_i^m\vert\vert\leq
4(1+||\tau_n||_{V^{\infty}})\va.$$

As $G$ is $\sigma$-compact there exist compact subsets $K_l$ such
that $K_l\subset K_{l+1}$ and $\cup_{l=1}^{\infty} K_l=G$, and
therefore for any $w\in T(G)$, $w=\lim w\chi_{K_l\times K_l}$ in
$T(G)$ giving that $\vert\vert\tau_nw-w\vert\vert\to 0$ for any
$w\in \Phi(E)$.
\end{proof}

\begin{cor}\label{union}
Any finite union  of sets of type (\ref{nontrian}) is a Ditkin
set.
\end{cor}

\begin{proof}
It follows from Proposition~\ref{nontr} and
\cite[Theorem~7.1]{sht}.\end{proof}

\begin{remark}
It follows, in particular,  from Corollary~\ref{union} that  finite
unions  of sets of type (\ref{nontrian}) are sets of operator
synthesis. This was also proved by Todorov in \cite{todorov} using
another method.
\end{remark}

 We will now establish a connection between (strong) Ditkin sets
for $A(G)$ and (strong) operator Ditkin sets.

 For $w\in
T(G)$, define as in \cite{var} $$Qw(s)=\int_Gw(sr,r)dr.$$ If $w
=\sum_{i=1}^{\infty}f_i\otimes g_i$ with
$\sum_{i=1}^{\infty}\vert\vert f_i\vert\vert_2^2\cdot
\sum_{i=1}^{\infty}\vert\vert g_i\vert\vert_2^2<\infty$. Then
$$Qw(s)=\sum_{i=1}^{\infty}(g_i*\check f_i)(s^{-1})\in A(G)$$ and,
moreover, $\vert\vert Qw\vert\vert\leq \vert\vert
w\vert\vert_{T(G)}$. Thus $Q:T(G)\to A(G)$ defines a contraction
operator.

\begin{theorem}\label{ditkin}
Let $G$ be a second countable locally compact group.
 If $E^*$ is a strong $m$-Ditkin set then $E$ is a local Ditkin set.
 If $E$ is a strong Ditkin set then $E^*$ is an $m$-Ditkin set.
\end{theorem}
\begin{proof}
Assume first that $E^*$ is a strong $m$-Ditkin set. Let
$\{\Psi_n\}_{n=1}^{\infty}\subset V^{\infty}(G)$ be a sequence
from the definition of a strong $m$-Ditkin set and let $u\in
I_c(E)$. For a compact subset $K\subset G$ containing the support
of $u$, define $(Nu)_K(s,t)=u(st^{-1})\chi_K(t)$. We have
$(Nu)_K(s,t)=\chi_{MK}(s)u(st^{-1})\chi_K(t)$, where $M=\supp
(u)$. As $u(st^{-1})\in V^{\infty}(G)$ and $|MK|<\infty$,
$|K|<\infty$, it yields $(Nu)_K\in T(G)$. Moreover, $(Nu)_K$
vanishes on $E^*$. Therefore,
$\vert\vert\Psi_n(Nu)_K-(Nu)_K\vert\vert_{T(G)}\to 0$ as
$n\to\infty$. Thus given $\varepsilon>0$, there exists $N$ such
that for $n>N$
$$\vert\vert
Q(\Psi_n(Nu)_K)-Q((Nu)_K)\vert\vert_A\leq
\vert\vert\Psi_n(Nu)_K-(Nu)_K\vert\vert_{T(G)}<\frac{\varepsilon}{\val
K}.$$

 On
the other hand,
\begin{eqnarray*}
Q(\Psi_n(Nu)_K)(s)&=&\int_G\Psi_n(sr,r)Nu(sr,r)\chi_K(r)dr
=\int_G\Psi_n(sr,r)u(s)\chi_K(r)dr\\&=&
u(s)\int_G\chi_{MK}(sr)\Psi_n(sr,r)\chi_K(r)dr=u(s)
Q(\Psi_n(\chi_{MK}\otimes\chi_K))(s)
\end{eqnarray*}
and similarly $Q((Nu)_K)=u\vert K\vert$ giving us
$$\vert\vert
Q(\Psi_n(Nu)_K)-Q((Nu)_K)\vert\vert_A=\vert\vert
uw_n-u|K|\vert\vert,$$

where $w_n=Q(\Psi_n(\chi_M\otimes\chi_K))$. We obtain now that for
$\tau_n=w_n/|K|$,
$$\vert\vert u\tau_n-u\vert\vert_A<\frac{\val K\varepsilon}{|K|}=\varepsilon.$$
Hence $E$ is a local Ditkin set.

Suppose $E$ is a strong Ditkin set. Let  $v_n\in J_A^0(G)$ be a
sequence such that $\vert\vert v_nf-f\vert\vert_A\to 0$ as
$n\to\infty$ for any $f\in I(E)$. By Proposition~\ref{nico},
$Nv_n\in V^{\infty}(G)$, $n=1,2,\ldots$, and
\begin{equation}\label{vinfty}
\vert\vert(Nv_n)(Nf)-Nf\vert\vert_{V^{\infty}}\leq\vert\vert v_nf-f\vert\vert_A\to
0.
\end{equation}
Next we show that $(Nv_n)w\cdot T\to w\cdot T$ ultra-weakly for
any compactly supported $T$ in $B(L_2(G))$ and any compactly
supported $w$ in $V^{\infty}(G)$ such that $w=0$ on $E^*$ using
arguments similar to the ones in the proof of
Theorem~\ref{ops}.

Assume $\supp (w)\subset K\times K$, for a
compact set $K\subset G$. Let $w_{kj}^{\pi}$ and $\tilde w_{kj}^{\pi}$
be $V^{\infty}$-functions as defined in (\ref{wkj}) and let
$u_{kj}^{\pi}$ be the matrix coefficient $(\pi(\cdot)e_k,e_j)$,
$\pi\in \hat G$. We have
$w_{kj}^{\pi}(s,t)=\sum_{l}u_{lj}^{\pi}(s)\tilde w_{kl}^{\pi}(s,t)$ by (\ref{wkjpi}).
Hence, for $u\in A(G)$, $\Psi\in T(G)$,
\begin{eqnarray*}
&&|\langle(Nv_n-1)w_{kj}^{\pi}\cdot T,N(u)\Psi\rangle|\leq\sum_{l}|\langle(Nv_n-1)\tilde
 w_{kl}^{\pi}\cdot T,u_{jl}^{\pi}N(u)\Psi\rangle|\\&&=\sum_{l\in L}|\langle(Nv_n-1)\tilde
 w_{kl}^{\pi}\cdot T,u_{jl}^{\pi}N(u)\Psi\rangle|+\sum_{l\notin L}|\langle(Nv_n-1)\tilde
 w_{kl}^{\pi}\cdot T,u_{jl}^{\pi}N(u)\Psi\rangle|
\end{eqnarray*}
where $L$ is a finite set. For the infinite sum we can apply the
estimate (\ref{inf}) with $(Nv_n-1)\Psi$ instead of $\Psi$,
 so that for given $\varepsilon>0$ there exists a
finite  subset $L$ such that
$$\sum_{l\notin L}|\langle(Nv_n-1)\tilde
 w_{kl}^{\pi}\cdot T,u_{jl}^{\pi}N(u)\Psi\rangle|=\sum_{l\notin L}|\langle\tilde
 w_{kl}^{\pi}\cdot T,u_{jl}^{\pi}N(u)(Nv_n-1)\Psi\rangle|<\epsilon$$
For the finite sum we note first that
$\tilde w_{kl}^{\pi}=N(v_{kl}^{\pi})$ for some $v_{kl}^{\pi}\in
M_{cb}A(G)$ by (\ref{equal}). Then  by
(\ref{vinfty}) there
exists $M>0$ such that for any $n>M$,
\begin{eqnarray*}
 \sum_{l\in L}|\langle(Nv_n-1)\tilde
 w_{kl}^{\pi}\cdot T,u_{jl}^{\pi}N(u)\Psi\rangle|=\sum_{l\in
 L}|\langle T,u_{jl}^{\pi}(Nv_n-1)N(v_{kl}^{\pi})N(u)\Psi\rangle|=\\
\sum_{l\in
 L}|\langle T,u_{jl}^{\pi}(Nv_n-1)N(v_{kl}^{\pi}u)\Psi\rangle|\leq\sum_{l\in
 L}||T||\cdot ||(Nv_n-1)N(v_{kl}^{\pi}u)||_{V^{\infty}}||\Psi||_{T(G)}
 <\varepsilon
\end{eqnarray*}
We obtain,
$$\langle(Nv_n-1)w_{kj}^{\pi}\cdot T,N(u)\Psi\rangle\to 0,\ n\to \infty.$$
As in the proof of Theorem~\ref{ops}, it yields
\begin{equation}\label{nvn}
\langle(Nv_n-1)w_{kj}^{\pi}\cdot T,\Psi\rangle\to 0,\ n\to \infty
\end{equation}
for any $\Psi\in T(G)$, since $\{Nv_n\}$ is bounded in norm.
Using the same arguments as in the end of the proof of
Theorem~\ref{synthesis} we have that for big enough compact set $K$
there exists a linear combination $u=\sum_{i=1}^dc_i u_i$ of matrix
coefficients such that
$$||w-u\chi_{K^{-1}K}\cdot w||_{T(G)}<\varepsilon.$$
Thus for $\omega\in T(G)$, $\supp(\omega)\subset K\times K$,
\begin{eqnarray*}
|\langle(Nv_n-1)w\cdot T,\omega\rangle|&\leq& |\langle(Nv_n-1)\cdot
T,\omega(w-u\chi_{K^{-1}K}\cdot w)\rangle|\\&+&|\langle(Nv_n-1)\cdot
T,\omega(u\chi_{K^{-1}K}\cdot w)\rangle|.
\end{eqnarray*}
The first summand is less than $C\varepsilon$ for some constant $C$,
for $n$ large enough ,  as $\{Nv_n\}$ is bounded in norm. As $\omega
(u_{kj}^{\pi}\cdot w)=\omega w_{kj}^{\pi}$ by (\ref{wkj}) for all
$k$, $j$, $\pi$,
 we have by (\ref{nvn})
$$\langle(Nv_n-1)\cdot
T,\omega(u_{kj}^{\pi}\chi_{K^{-1}K}\cdot w)\rangle=\langle(Nv_n-1)w_{kj}^{\pi}\cdot
T,\omega\rangle\to 0, n\to \infty,$$
and hence
$|\langle(Nv_n-1)\cdot
T,\omega(u\chi_{K^{-1}K}\cdot w)\rangle|<\varepsilon$
for $n$ large enough.
Thus, $$\langle(Nv_n-1)w\cdot T,\omega\rangle\to 0$$ for each $\omega\in
T(G)$.
 Hence
$(Nv_n)w\omega\to w\omega$ weakly for each $\omega\in T(G)$. As
the set of all linear combinations of $w\omega$, where $w\in
V^{\infty}$, $w=0$ on $E^*$ and $\omega \in T(G)$, is dense in
$\Phi(E)$ (\cite[Proposition~5.3]{sht2}), $(Nv_n)\omega\to \omega$
weakly for any $\omega\in\Phi(E)$.  Taking if necessary convex
linear combinations of the $Nv_n$'s,  we obtain elements $w_n \in
V^\infty(G), n\in \N,$ such that $\vert\vert
w_n\omega-\omega\vert\vert_{T(G)}\to 0$ as $n$ tends to infinity.
Clearly, the sequence $\{w_n\}_n$ satisfies the necessary
conditions.
 \end{proof}

\begin{cor}
Any closed subgroup $H$ of $G$ is a local Ditkin set.
\end{cor}
\begin{proof}
We have $H^*=\{(s,t):f(s)=f(t)\}$, where $f:G\to G\backslash H$,
$t\mapsto tH$. As $G$ is metrizable, by
\cite[((8.14)]{hewitt-rossI}, $G\backslash H$ is metrizable. The
statement now follows from Proposition~\ref{nontr} and
Theorem~\ref{ditkin}.
\end{proof}

For $G$ amenable the statement was obtained in \cite{fkls} showing
that $I_A(H)$ has a bounded approximate identity. For arbitrary
locally compact $G$ and $H$ a closed neutral subgroup it was proved
in \cite{derighetti}. Note  that there are
   closed subgroups of l.c.s.c. groups which are not neutral
(see \cite[p.107]{roelcke}).

\vspace{1cm}

{\bf Acknowledgments.}

\vspace{0.5cm}

 We are indebted to A.T.Lau, V.Shulman and  N.Spronk
for helpful discussions and valuable information.

\end{document}